\documentclass{article}

\usepackage{graphicx}
\usepackage[fleqn]{amsmath}
\usepackage{amsthm}
\usepackage{amssymb}
\usepackage{amsbsy}
\usepackage{amsfonts}
\usepackage{amstext}
\usepackage{amscd}
\usepackage[dvips]{epsfig}
\allowdisplaybreaks

\ifx\undefined\bysame
\newcommand{\bysame}{\leavevmode\hbox to3em{\hrulefill}\,}
\fi

\newcommand{\cA}{{\cal A}}
\newcommand{\e}{\varepsilon}
\newcommand{\hT}{{\hat\Theta}}
\newcommand{\Q}{{\Bbb Q}}
\newcommand{\Z}{{\Bbb Z}}

\newcommand{\pc}[1]{\mbox{$\begin{array}{c}
   \includegraphics[scale=0.5]{figure/#1.ps}
   \end{array}$}}

\theoremstyle{plain}
   \newtheorem{thm}{Theorem}[section]
   \newtheorem{lem}[thm]{Lemma}
   \newtheorem{prop}[thm]{Proposition}
   \newtheorem{cor}[thm]{Corollary}
\theoremstyle{definition}
   
   \newtheorem{ex}[thm]{Example}
   \newtheorem{rem}[thm]{Remark}

\begin{document}

\vspace*{2pc}

\Large
\centerline{\bf\LARGE A cabling formula for the 2-loop polynomial of knots}

\vskip 2pc

\centerline{Tomotada Ohtsuki}

\begin{abstract}
The 2-loop polynomial is a polynomial presenting the 2-loop part
of the Kontsevich invariant of knots.
We show a cabling formula for the 2-loop polynomial of knots.
In particular, we calculate the 2-loop polynomial for torus knots.
\end{abstract}

\bigskip
\large

\baselineskip 1.4pc

The Kontsevich invariant is a very strong invariant of knots
(which dominates all quantum invariants and all Vassiliev invariants)
and it is expected that the Kontsevich invariant will classify knots.
A problem when we study the Kontsevich invariant is that
it is difficult to calculated the Kontsevich invariant of 
an arbitrarily given knot concretely.
It has recently been shown 
\cite{Rozansky_rc,Kricker_rc,GK_rat}\footnote{\small
It was conjectured by Rozansky \cite{Rozansky_rc}.
The existence of such rational presentations has been proved by 
Kricker \cite{Kricker_rc}
(though such a rational presentation itself is not necessarily
a knot invariant in a general loop degree).
Further, Garoufalidis and Kricker \cite{GK_rat}
defined a knot invariant in any loop degree, 
from which such a rational presentation can be deduced.}
that the infinite sum of the terms of the Kontsevich invariant
with a fixed loop number is presented by using polynomials
(after appropriate normalization by the Alexander polynomial).
In particular, it is known\footnote{\small
This follows from the theory of \cite{BG_MMR} on the MMR conjecture.
See also \cite{Kricker_rc,GK_rat} and references therein.}
that the 1-loop part is presented by the Alexander polynomial.
The polynomial giving the 2-loop part is called the 2-loop polynomial.
The values of the 2-loop polynomial 
has been calculated so far only for particular\footnote{\small
A table of the 2-loop polynomial for knots with up to 7 crossings
is given by Rozansky \cite{Rozansky_2lp}.
The 2-loop polynomial of knots with the trivial Alexander polynomial
can often been calculated by surgery formulas \cite{GK_rat,Kricker_s}.} 
classes of knots.

\begin{figure}[htpb]
$$
\pc{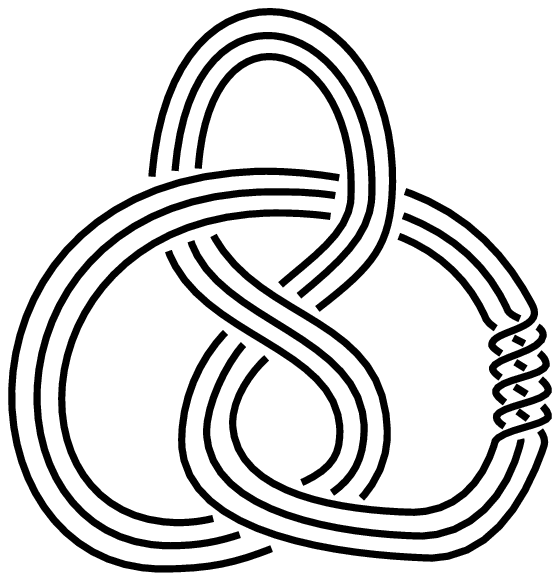} \mbox{ is the $(3,5)$ cable knot of }
\pc{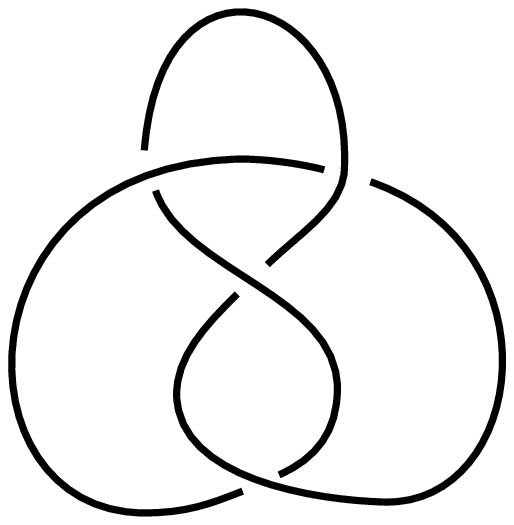} .
$$
\caption{\label{fig.cable}
A cable knot of a knot}
\end{figure}

In this paper,
we give a cabling formula for the 2-loop polynomial
(Theorem \ref{thm.cabling}),
which presents the 2-loop polynomial of a cable knot 
(see Figure \ref{fig.cable})
of a knot $K$
in terms of the 2-loop polynomial of $K$.
In particular, we calculate a formula of the 2-loop polynomial for torus knots
(Theorem \ref{thm.2ltorus}).
This formula and the cabling formula
are also obtained independently by March\'e \cite{Marche,Marche_PhD}.

This paper is organized as follows.
In Section 1
we review the definition of the 2-loop polynomial.
In Section 2
we calculate the 2-loop polynomial of torus knots
as the 2-loop part of the cabling formula
of the Kontsevich invariant of the trivial knot.
In Section 3
we give a cabling formula for the 2-loop polynomial.
In Section 4
we show relations to some Vassiliev invariants.

The author would like to thank
Andrew Kricker, Thang Le, Lev Rozansky, Julien March\'e,
Stavros Garoufalidis, Dror Bar-Natan
for valuable discussions and comments.

\section{The Kontsevich invariant and the 2-loop polynomial}

The 2-loop polynomial is a polynomial presenting
the 2-loop part of the Kontsevich invariant.
In this section,
we review its definition and a cabling formula of the Kontsevich invariant.

\begin{figure}[htpb]
$$
\begin{array}{lcl}
\mbox{The AS relation:} &&
\pc{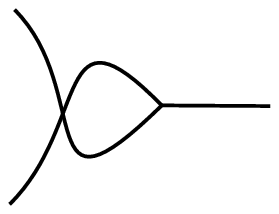} = - \pc{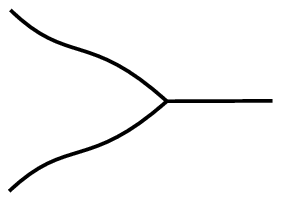} \\
\mbox{The IHX relation:} &&
\pc{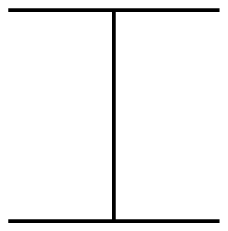} = \pc{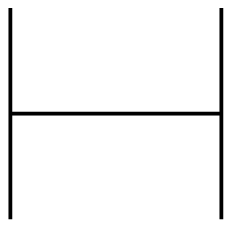} - \pc{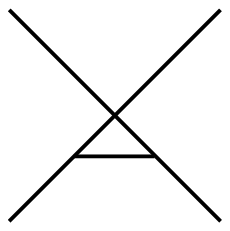}
\end{array}
$$
\caption{\label{fig.as-ihx}
The AS and IHX relations}
\end{figure}

An {\it open Jacobi diagram} is a uni-trivalent graph
such that a cyclic order of the three edges
around each trivalent vertex of the graph is fixed.
Let $\cA(\ast)$ be the vector space over $\Q$
spanned by open Jacobi diagrams
subject to the AS and IHX relations;
see Figure \ref{fig.as-ihx} for the relations.

The {\it Kontsevich invariant} $Z^\sigma(K)$ of a framed knot $K$
is defined in $\cA(\ast)$; for a definition\footnote{\small
In literatures, the Kontsevich invariant is
often defined by $Z(K)$ in the space $\cA(S^1)$.
The version $Z^\sigma(K)$ is defined to be the image of $Z(K)$
by the inverse map $\sigma$ of the Poincare-Birkhoff-Witt isomorphism
$\cA(\ast) \to \cA(S^1)$.}
see {\it e.g.} \cite{Ohtsuki_book}.
It is known \cite{LMO}
that the value of the Kontsevich invariant for each knot is group-like,
which implies that
it is presented by the exponential of some primitive element.
That is,
$Z^\sigma(K)$ is presented by the exponential of a primitive element,
where a {\it primitive element} of $\cA(\ast)$
is a linear sum of connected open Jacobi diagrams.

For example, it is shown \cite{BLT} that
the Kontsevich invariant of the trivial knot, denoted by $\Omega$, 
is presented by
$$
Z^\sigma(\mbox{the trivial knot}) = \Omega = \exp_\sqcup(\omega),
$$
where $\exp_\sqcup$ denotes the exponential with respect to
the disjoint-union product, and $\omega$ is defined by
$$
\omega = 
{\begin{array}{c}
{\begin{picture}(105,64)
\put(21,53){$\frac12 \log \frac{\sinh(x/2)}{x/2}$}
\put(0,20){\pc{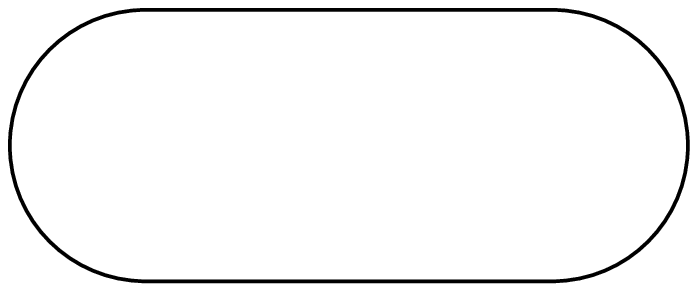}}
\end{picture}} \end{array}} .
$$
Here, a label of a power seires
$f(x)=c_0+c_1x+c_2 x^2+c_3x^3+\cdots$
implies
$$
\pc{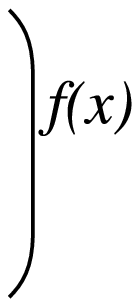}= c_0 \pc{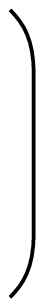}+c_1 \pc{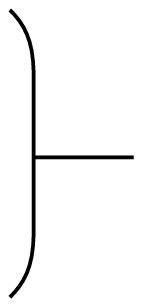}
+c_2 \pc{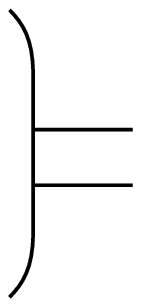}+c_3 \pc{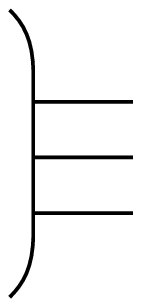}+\cdots,
$$
where a label is put on either of the sides of an edge,
and the corresponding legs are written 
in the same side of the edge.\footnote{\small
Our notation is different from the notation in \cite{GK_rat,Kricker_s}
where a label of an edge is defined by setting 
a local orientation of the edge
that determines the side in which we write the corresponding legs.}
Note that
$\pc{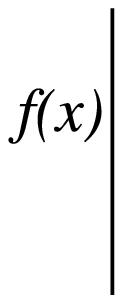}=\pc{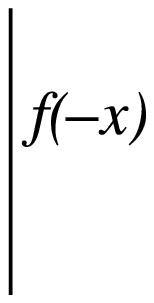}$
by the AS relation,
in the notation of this paper.

Let $K$ be a framed knot with $0$ framing.
(Throughout this paper, we often mean 
a framed knot with $0$ framing also by a knot, abusing terminology.)
The {\it loop expansion} of the Kontsevich invariant is given by
\begin{align*}
\log_\sqcup Z^\sigma(K) &= {\begin{array}{c}
{\begin{picture}(155,64)
\put(12,53){$\frac12 \log \frac{\sinh(x/2)}{x/2}
-\frac12 \log \Delta{\!}_K{\!}(e^x)$}
\put(0,20){\pc{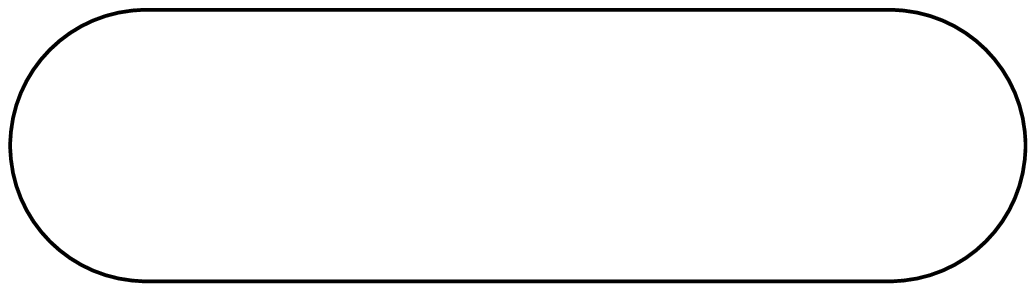}}
\end{picture}} \end{array}} 
+ \sum_i^{\mbox{\scriptsize finite}} \!\!\! {\begin{array}{c}
{\begin{picture}(120,80)
\put(27,71){$p_{i,1}(e^x)/\Delta{\!}_K{\!}(e^x)$}
\put(27,39){$p_{i,2}(e^x)/\Delta{\!}_K{\!}(e^x)$}
\put(27,7){$p_{i,3}(e^x)/\Delta{\!}_K{\!}(e^x)$}
\put(0,28){\pc{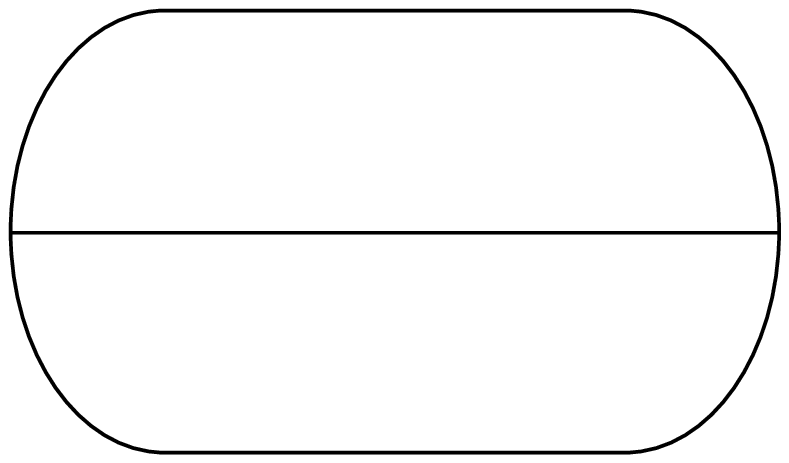}}
\end{picture}} \end{array}}
\\
&+ \ 
\big( \mbox{terms of $(\ge 3)$-loop} \big),
\end{align*}
where 
$\log_\sqcup$ denotes the logarithm
with respect to the disjoint-union product,
and $\Delta_K(t)$ is the normalized\footnote{\small
We suppose that $\Delta_K(t)$ is normalized, satisfying that
$\Delta_K(t)=\Delta_K(t^{-1})$ and $\Delta_K(1)=1$.}
Alexander polynomial of $K$,
and $p_{i,j}(e^x)$ is a polynomial in $e^x$.
The 2-loop part is characterized by the polynomial,
$$
\Theta'_K(t_1,t_2,t_3) = \sum_i p_{i,1}(t_1) p_{i,2}(t_2) p_{i,3}(t_3).
$$
We call its symmetrization,\footnote{\small
With respect to the symmetry of the theta graph, of order 12.}
$$
\Theta_K(t_1,t_2,t_3) = 
\!\!\!\!\!\sum_{\substack{\e=\pm1 \\ \{i,j,k \}=\{1,2,3\} }}\!\!\!\!\!
\Theta'_K(t^\e_{i},t^\e_{j},t^\e_{k})
\quad
\in \Q[t_1^{\pm1},t_2^{\pm1},t_3^{\pm1}]/(t_1 t_2 t_3 =1),
$$
the {\it 2-loop polynomial} of $K$,
which is an invariant\footnote{\small
This is not trivial,
since there is another 2-loop trivalent graph,
what is called, a ``dumbbell diagram''.}
of $K$.
(Note that this normalization of $\Theta_K(t_1,t_2,t_3)$ is
12 times the usual normalization.)
$\Theta_K(t,t^{-1},1)$ is a symmetric polynomial in $t^{\pm1}$
divisible by $t-1$ (since $\Theta_K(1,1,1)=0$)
and, hence, divisible by $(t-1)^2$.
We define the {\it reduced 2-loop polynomial} by
$$
\hT_K(t)
=\frac{\Theta_K(t,t^{-1},1)}{(t^{1/2}-t^{-1/2})^2}
\quad \in \Q[t^{\pm1}],
$$
which is a symmetric polynomial in $t^{\pm1}$.

Let us review the cabling formula of the Kontsevich invariant
of \cite{BLT}.
Another version of the Kontsevich invariant,
called the {\it wheeled Kontsevich invariant} \cite{BarNatan-Lawrence},
is defined by
$$
Z^w(K) = \partial_\Omega^{-1} Z^\sigma(K),
$$
where $\partial_\Omega:\cA(\ast) \to \cA(\ast)$
is the {\it wheeling isomorphism};
see \cite{BLT}.
Here, for open Jacobi diagrams $C$ and $D$,
$\partial_C(D)$ is defined to be $0$
if $C$ has more univalent vertices than $D$,
and the sum of all ways of gluing all univalent vertices of $C$
to some univalent vertices of $D$ otherwise.
We graphically present it by
$$
\partial_C(D) = \pc{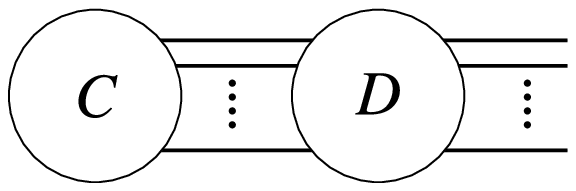}.
$$
Let $\Psi^{(p)}: \cA(\ast) \to \cA(\ast)$
be the map which takes a diagram with $k$ univalent vertices
to its $p^k$ multiple.
The {\it $(p,q)$ cable knot} of a knot $K$
is the knot given by a simple closed curve
on the boundary torus of a tubular neighborhood of $K$
which winds $q$ times in the meridian direction
and $p$ times in the longitude direction
(see {\it e.g.} \cite{Lickorish_book});
for example see Figure \ref{fig.cable}.
The cabling formula of the Kontsevich invariant is given by\footnote{\small
Proposition \ref{prop.cabling} is obtained from
Theorem 1 of \cite{BLT} by pulling back by the isomorphism
$\cA(\ast) \overset{\partial_\Omega}\longrightarrow
\cA(\ast) \overset{\chi}\longrightarrow \cA(S^1)$,
and by modifying the contribution from the framing of the cable knot,
noting that the $(p,q)$ cable knot in the definition of \cite{BLT}
has framing $(p-1)q$.}

\begin{prop}[Le {\rm (\cite{BLT}, see also \cite{Willerton})}]
\label{prop.cabling}
Let $K$ be a framed knot with $0$ framing,
and let $K^{(p,q)}$ be the $(p,q)$ cable knot of $K$ (with $0$ framing).
Then, 
$$
Z^w(K^{(p,q)}) = \partial_\Omega^{-1} \Psi^{(p)} \partial_\Omega \Big(
Z^w(K) \sqcup
\exp_\sqcup \big( \frac{q}{2p} \pc{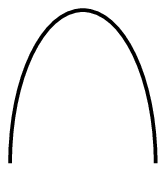} -\frac{q}{48p} \theta 
\big) \Big)
\sqcup
\exp_\sqcup \big( -\frac{p q}{2} \pc{arc} +\frac{p q}{48} \theta \big) .
$$
\end{prop}

\section{The 2-loop polynomial of a torus knot}

In this section,
we calculate the 2-loop polynomial of a torus knot,
picking up the 2-loop part
of the cabling formula of the Kontsevich invariant of the trivial knot.
The 2-loop part of the Kontsevich invariant for torus knots
is also calculated\footnote{\small
Bar-Natan has also obtained some presentation of 
the wheeled Kontsevich invariant for torus knots (private communication).}
independently by March\'e \cite{Marche,Marche_PhD}.

\begin{figure}[htpb]
$$
\begin{picture}(100,50)
\put(10,10){\pc{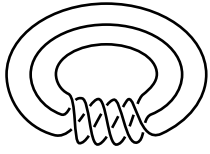}}
\end{picture}
$$
\caption{\label{fig.torus_knot}
The $(5,3)$ torus knot}
\end{figure}

The {\it torus knot} $T(p,q)$ of type $(p,q)$
is the $(p,q)$ cable knot of the trivial knot
(which is isotopic to $T(q,p)$);
for example see Figure \ref{fig.torus_knot}.
It is known, see {\it e.g.} \cite{Lickorish_book},
that the Alexander polynomial of a torus knot is given by
$$
\Delta_{T(p,q)}(t) =
\frac{(t^{p q/2}-t^{-p q/2})(t^{1/2}-t^{-1/2})}
     {(t^{p/2}-t^{-p/2})(t^{q/2}-t^{-q/2})}.
$$

\begin{thm}
\label{thm.2ltorus}
The 2-loop polynomial of the torus knot $T(p,q)$ of type $(p,q)$
is given by\footnote{\small
This value coincides with the value in \cite{Marche,Marche_PhD}.
However, the values of the 2-loop polynomial
for some torus knots in Table 2 of \cite{Rozansky_2lp}
have opposite signs to our values.
The signs of some values in Table 2 of \cite{Rozansky_2lp}
might not be correct.}
$$
\Theta_{T(p,q)}(t_1,t_2,t_3)
= - \frac1{4} \sum_{ \{i,j,k\} = \{1,2,3 \} } \!\!\!\!\!
\psi_{p,q}(t_i) \psi_{q,p}(t_j) \Delta_{T(p,q)}(t_k)
\quad \in \Z[t_1^{\pm1},t_2^{\pm1},t_3^{\pm1}]/(t_1t_2t_3=1),
$$
where $\psi_{p,q}$ is defined by
\begin{align*}
\psi_{p,q}(t) &= \Delta_{T(p,q)}(t) \cdot
\Big( \frac{t^{p/2}+t^{-p/2}}{t^{p/2}-t^{-p/2}} 
-q \cdot \frac{t^{p q/2}+t^{-p q/2}}{t^{p q/2}-t^{-p q/2}} \Big) \\*
&= 
\frac{t^{1/2}-t^{-1/2}}{(t^{p/2}-t^{-p/2})(t^{q/2}-t^{-q/2})}
\Big(
(t^{p/2}+t^{-p/2})\cdot\frac{t^{p q/2}-t^{-p q/2}}{t^{p/2}-t^{-p/2}} 
-q (t^{p q/2}+t^{-p q/2})
\Big).
\end{align*}
In particular,
$\Theta_{T(p,q)}(t_1,t_2,t_3)$ is a polynomial in 
$t_1^{\pm1}$, $t_2^{\pm1}$, $t_3^{\pm1}$
with integer coefficients of
$\mbox{degree}_{t_1} \big( \Theta_{T(p,q)}(t_1,t_2,t_1^{-1}t_2^{-1}) \big)
= (p-1)(q-1)$.
\end{thm}

\begin{rem}
$\psi_{p,q}(t)$ is not a polynomial, but a rational function,
while $\Theta_{T(p,q)}(t_1,t_2,t_3)$ is a polynomial.
Rozansky \cite{Rozansky_2lp} suggests that
the 2-loop polynomial is a polynomial with integer coefficients;
this holds for torus knots by the theorem.
He also suggests a conjectural inequality
$$
\mbox{degree}_{t_1} \big( \Theta_{K}(t_1,t_2,t_1^{-1}t_2^{-1}) \big)
\le 2 g(K),
$$
where $g(K)$ denotes the genus of $K$.
Since the genus of $T(p,q)$ equals
$(p-1)(q-1)/2$ (see {\it e.g.} \cite{Lickorish_book}),
torus knots give the equality of the above formula.
\end{rem}

\begin{rem}
The $sl_2$ reduction of the $n$-loop part of the Kontsevich invariant
of the Kontsevich invariant is equal to the $n$th line
in the expansion of the colored Jones polynomial.
Rozansky \cite{Rozansky_hot} has calculated it for torus knots.
\end{rem}

For $\alpha, \beta \in \cA(\ast)$
we write $\alpha\equiv\beta$
if $\alpha-\beta$ is equal to a linear sum of Jacobi diagrams,
either, of $(\ge 3)$-loop, or,
having a component of a trivalent graph
({\it i.e.}, a component with no univalent vertices).

\begin{proof}[Proof of Theorem \ref{thm.2ltorus}]
Since the torus knot $T(p,q)$ is obtained from the trivial knot by cabling,
we have that
$$
Z^w\big( T(p,q) \big)
\equiv
\partial_\Omega^{-1} \Psi^{(p)} \partial_\Omega
\Big( \Omega \sqcup
\exp_\sqcup \big( \frac{q}{2p} \pc{arc} \big) \Big)
\sqcup 
\exp_\sqcup \big( -\frac{p q}{2} \pc{arc} \big) 
$$
by Proposition \ref{prop.cabling}.
The first term of the right hand side is calculated as follows.
From the definition of $\partial_\Omega$,
\begin{equation}
\label{eq.O-O-e}
\partial_\Omega \Big( 
\exp_\sqcup \big( \frac{q}{2p} \pc{arc} \big) \sqcup \Omega 
\Big)
=
\pc{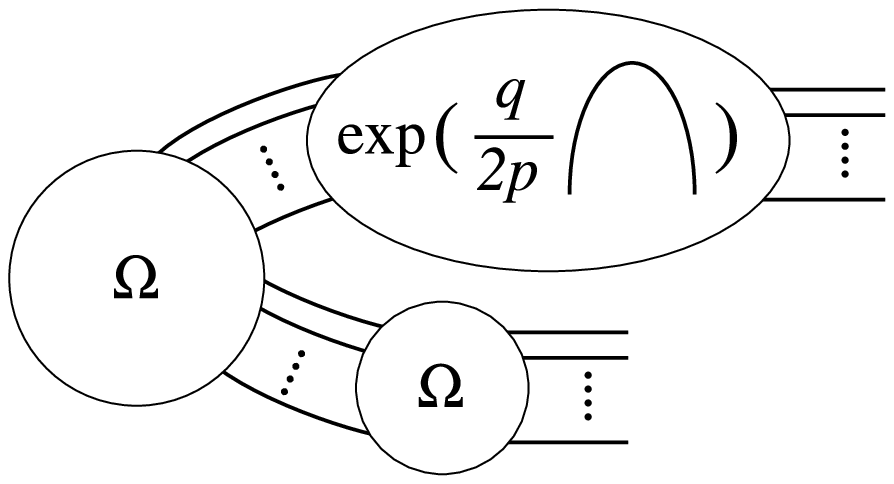}.
\end{equation}
Since any component of $\Omega$ has a loop,
the $(\le 1)$-loop part of the right hand side has no edges 
between the two $\Omega$'s,
and, hence, this part is presented by
$$
\partial_\Omega \exp_\sqcup \big( \frac{q}{2p} \pc{arc} \big)
\sqcup \Omega .
$$
Further, its first term is given by
$$
\partial_\Omega \exp_\sqcup \big( \frac{q}{2p} \pc{arc} \big)
\equiv
\exp_\sqcup \big( \frac{q}{2p} \pc{arc} \big)
\sqcup \Omega_{\frac{q}{p} x} ,
$$
where the equivalence is obtained in the same was as
Lemma 6.3 of \cite{BLT}.
The primitive part of the 2-loop part of the right hand side of
(\ref{eq.O-O-e})
is equal to a linear sum of diagrams, each of which
has precisely one edge between the two $\Omega$'s.
Hence, it is presented by
$$
\pc{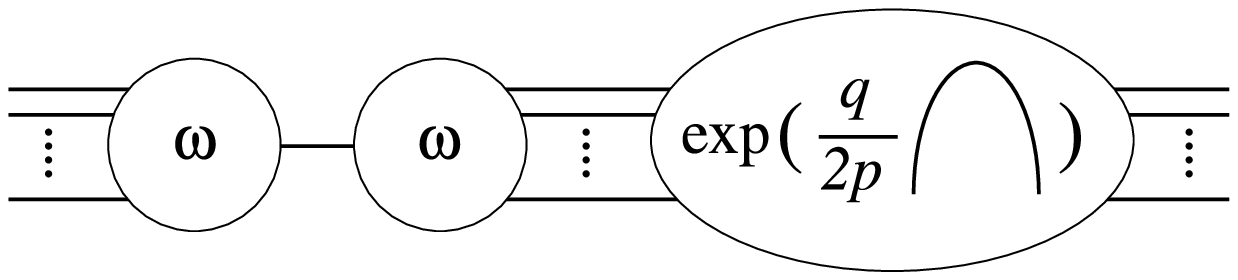}.
$$
Since
$$
\pc{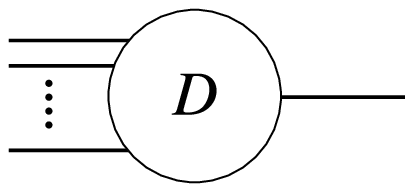} = \pc{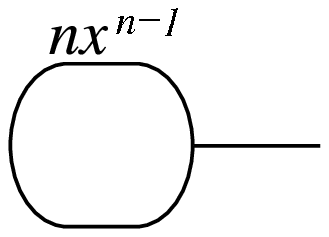}
\qquad \mbox{ for }
D = \pc{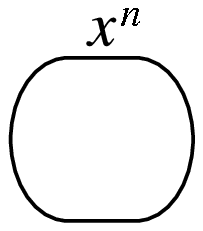},
$$
the previous diagram is equivalent to
$$
\begin{array}{c}\begin{picture}(90,30)
\put(0,5){\pc{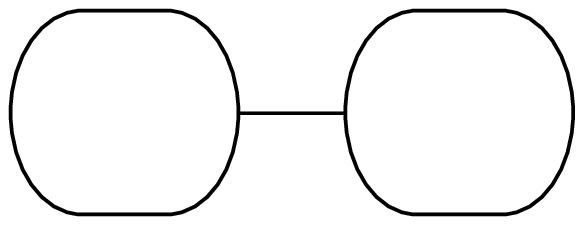}}
\put(11,30){$f(x)$}
\put(55,32){$f( \frac{q}{p} x)$}
\end{picture}\end{array},
$$
where $f(x)$ is given by
$$
f(x) = \frac{d}{d x} \Big( \frac12 \log \frac{\sinh x/2}{x/2} \Big)
= \frac14 \cdot \frac{e^{x/2}+e^{-x/2}}{e^{x/2}-e^{-x/2}} - \frac1{2x}.
$$
Hence, the $(\le2)$-loop part of (\ref{eq.O-O-e}) is presented by
\begin{equation}
\label{eq.dOeO}
\partial_\Omega \Big( 
\exp_\sqcup \big( \frac{q}{2p} \pc{arc} \big) \sqcup \Omega 
\Big)
\equiv
\exp \big( \frac{q}{2p} \pc{arc} \big)
\sqcup \Omega \sqcup \Omega_{\frac{q}{p} x} \sqcup 
\exp_\sqcup \Big( \! \begin{array}{c}\begin{picture}(90,30)
\put(0,5){\pc{megane3}}
\put(11,30){$f(x)$}
\put(55,32){$f( \frac{q}{p} x)$}
\end{picture}\end{array} \!
\Big).
\end{equation}
The map $\Psi^{(p)}$ sends this to
$$
\exp \big( \frac{p q}{2} \pc{arc} \big)
\sqcup \Omega_{p x} \sqcup \Omega_{q x} \sqcup 
\exp_\sqcup \Big( \! \begin{array}{c}\begin{picture}(90,30)
\put(0,5){\pc{megane3}}
\put(11,30){$f(p x)$}
\put(55,30){$f(q x)$}
\end{picture}\end{array} \!
\Big).
$$
Further, $\partial_\Omega^{-1}$ sends this (modulo the equivalence) to
$$
\partial_{\Omega^{-1}} \Big(
\exp \big( \frac{p q}{2} \pc{arc} \big)
\sqcup \Omega_{p x} \sqcup \Omega_{q x} \Big)
\sqcup 
\exp_\sqcup \Big( \! \begin{array}{c}\begin{picture}(90,30)
\put(0,5){\pc{megane3}}
\put(11,30){$f(p x)$}
\put(55,30){$f(q x)$}
\end{picture}\end{array} \!
\Big).
$$
Its first term is graphically shown as
\begin{equation}
\label{eq.O-OOe}
\pc{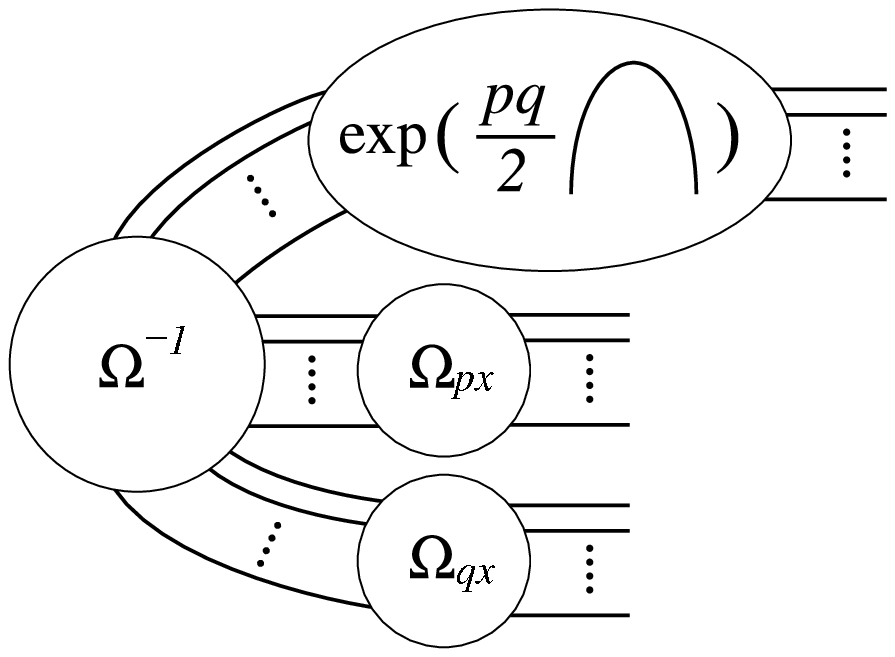}.
\end{equation}
The primitive part of the 2-loop part of this diagram is calculated
similarly as before;
for example,
when there is precisely one edge between $\Omega^{-1}$ and $\Omega_{p x}$,
we have the following component,
$$
\pc{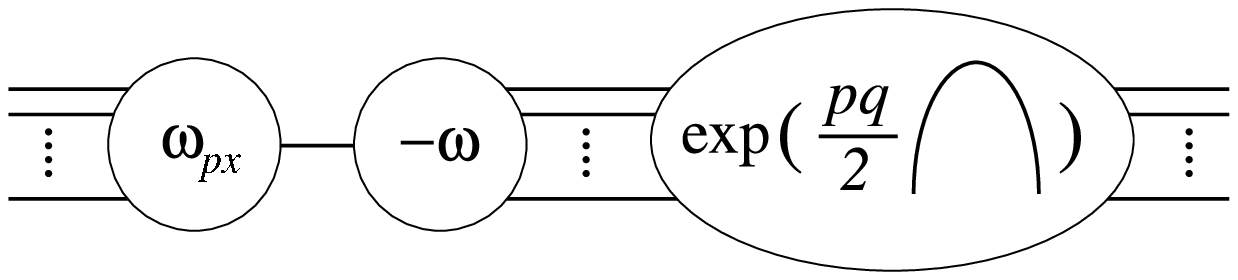}
\equiv
-p \!\!
\begin{array}{c}\begin{picture}(90,30)
\put(0,5){\pc{megane3}}
\put(8,30){$f(p x)$}
\put(52,30){$f(p q x)$}
\end{picture}\end{array}.
$$
Thus, the primitive part of the 2-loop part of (\ref{eq.O-OOe})
is equal to
\begin{align*}
& \Big(\mbox{the primitive part of the 2-loop part of 
$\partial_{\Omega^{-1}} 
\exp \big( \frac{p q}{2} \pc{arc} \big)$}\Big) \\
&\qquad
-p \!\!
\begin{array}{c}\begin{picture}(90,30)
\put(0,5){\pc{megane3}}
\put(8,30){$f(p x)$}
\put(52,30){$f(p q x)$}
\end{picture}\end{array}
-q \!\!
\begin{array}{c}\begin{picture}(90,30)
\put(0,5){\pc{megane3}}
\put(8,30){$f(q x)$}
\put(52,30){$f(p q x)$}
\end{picture}\end{array} \\
&=p q \!\!
\begin{array}{c}\begin{picture}(90,38)
\put(0,5){\pc{megane3}}
\put(3,30){$f(p q x)$}
\put(52,30){$f(p q x)$}
\end{picture}\end{array}
-p \!\!
\begin{array}{c}\begin{picture}(90,38)
\put(0,5){\pc{megane3}}
\put(8,30){$f(p x)$}
\put(52,30){$f(p q x)$}
\end{picture}\end{array}
-q \!\!
\begin{array}{c}\begin{picture}(90,38)
\put(0,5){\pc{megane3}}
\put(8,30){$f(q x)$}
\put(52,30){$f(p q x)$}
\end{picture}\end{array},
\end{align*}
where the equality is obtained from Lemma \ref{lem.pO-} below.
Hence, the primitive part of the 2-loop part of 
$Z^w \big( T(p,q) \big)$ is given by
\begin{align}
\!\!\!\!\!\!\!\!
& \begin{array}{c}\begin{picture}(90,35)
\put(0,5){\pc{megane3}}
\put(8,30){$f(p x)$}
\put(57,30){$f(q x)$}
\end{picture}\end{array}
+p q \!\!
\begin{array}{c}\begin{picture}(90,35)
\put(0,5){\pc{megane3}}
\put(3,30){$f(p q x)$}
\put(52,30){$f(p q x)$}
\end{picture}\end{array}
-p \!\!
\begin{array}{c}\begin{picture}(90,35)
\put(0,5){\pc{megane3}}
\put(8,30){$f(p x)$}
\put(52,30){$f(p q x)$}
\end{picture}\end{array}
-q \!\!
\begin{array}{c}\begin{picture}(90,35)
\put(0,5){\pc{megane3}}
\put(8,30){$f(q x)$}
\put(52,30){$f(p q x)$}
\end{picture}\end{array} \notag \\
\label{eq.me-theta}
\!\!\!\!\!\!\!\!&=
\frac1{16}
\begin{array}{c}
\begin{picture}(110,30)
\put(0,5){\pc{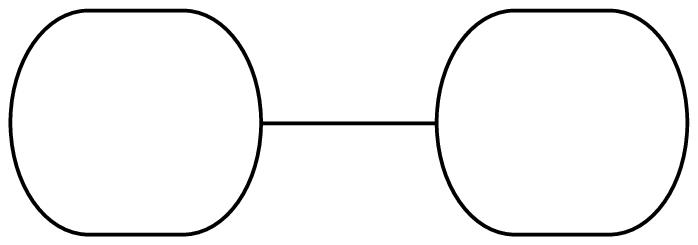}}
\put(10,32){$\phi_{p,q}(t)$}
\put(70,32){$\phi_{q,p}(t)$}
\end{picture}\end{array}
= - \frac18
\begin{array}{c}
\begin{picture}(80,65)
\put(0,20){\pc{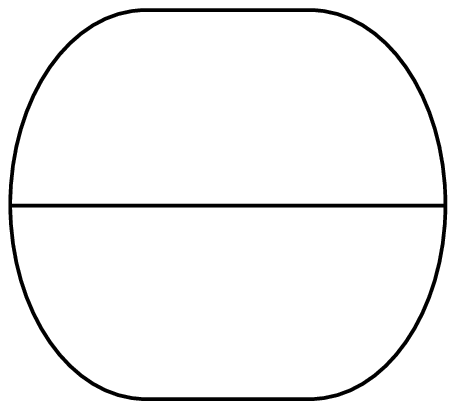}}
\put(22,60){$\phi_{p,q}(t)$}
\put(22,32){$\phi_{q,p}(t)$}
\end{picture}\end{array},
\end{align}
where we put $t=e^x$ and
$\phi_{p,q}$ is defined by
$\phi_{p,q}(e^x) = 4 \big( f(p x)-q f(p q x) \big)$, that is,
$$
\phi_{p,q}(t) = \frac{t^{p/2}+t^{-p/2}}{t^{p/2}-t^{-p/2}} 
-q \cdot \frac{t^{p q/2}+t^{-p q/2}}{t^{p q/2}-t^{-p q/2}}.
$$
Therefore, from the definition of the 2-loop polynomial,
we obtain the required formula.

By Corollary \ref{cor.torus} below,
the degree of $\hat\Theta_{T(p,q)}(t)$ equals $(p-1)(q-1)-1$.
Since 
$(t^{1/2}-t^{-1/2})^2 \hat\Theta_{T(p,q)}(t)=\Theta_{T(p,q)}(t,1,t^{-1})$ 
by definition,
$\mbox{degree}_{t_1} \big( \Theta_{T(p,q)}(t_1,t_2,t_1^{-1}t_2^{-1}) \big)$
is at least $(p-1)(q-1)$.
We can show that it is exactly $(p-1)(q-1)$
in the same way as the proof of Example \ref{ex.p2torus}.
\end{proof}

\begin{cor}
\label{cor.torus}
The reduced 2-loop polynomial of the torus knot $T(p,q)$ is given by
\begin{align*}
\hat\Theta_{T(p,q)}(t)
&= \frac1{2 (t^{1/2}-t^{-1/2})^2}   \psi_{p,q}(t) \psi_{q,p}(t) \\*
&= \frac12 \cdot \frac1{(t^{p/2}-t^{-p/2})^2} \cdot
\Big(
(t^{p/2}+t^{-p/2})\cdot\frac{t^{p q/2}-t^{-p q/2}}{t^{p/2}-t^{-p/2}} 
-q (t^{p q/2}+t^{-p q/2})
\Big) \\*
&\ \quad\times 
\frac1{(t^{q/2}-t^{-q/2})^2} \cdot
\Big(
(t^{q/2}+t^{-q/2})\cdot\frac{t^{p q/2}-t^{-p q/2}}{t^{q/2}-t^{-q/2}} 
-p (t^{p q/2}+t^{-p q/2})
\Big) .
\end{align*}
\end{cor}

\begin{lem}
\label{lem.pO-}
For a scalar $c$,
$$
\partial_\Omega^{-1} \exp \Big( \frac{c}2 \pc{arc} \Big)
\equiv
\exp \Big( \frac{c}2 \pc{arc} \Big)
\sqcup \Omega_{c x}^{-1}
\sqcup \exp_\sqcup \Big( c \!\!
\begin{array}{c}\begin{picture}(90,30)
\put(0,5){\pc{megane3}}
\put(8,30){$f(c x)$}
\put(57,30){$f(c x)$}
\end{picture}\end{array} \!
\Big).
$$
\end{lem}

\begin{proof}
From the definition of $\partial_\Omega$,
\begin{equation}
\label{eq.OOec}
\partial_\Omega \Big( 
\exp \big( \frac{c}2 \pc{arc} \big)
\sqcup \Omega_{c x}^{-1} \Big)
=
\pc{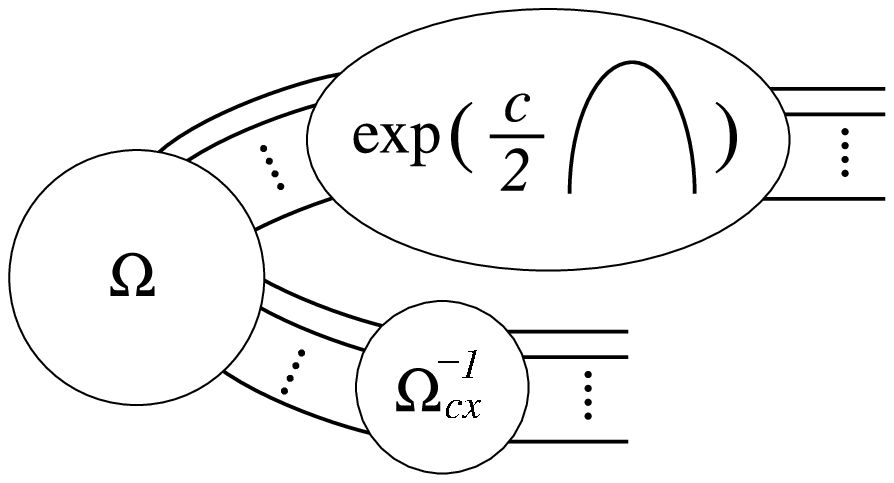}.
\end{equation}
Similarly as in the proof of Theorem \ref{thm.2ltorus},
the $(\le 1)$-loop part of the right hand side is presented by
$$
\partial_\Omega 
\exp \Big( \frac{c}2 \pc{arc} \Big)
\sqcup \Omega_{c x}^{-1} 
\equiv
\exp \Big( \frac{c}2 \pc{arc} \Big) .
$$
Further, the primitive part of the 2-loop part of the right hand side of
(\ref{eq.OOec}) is presented by
$$
\pc{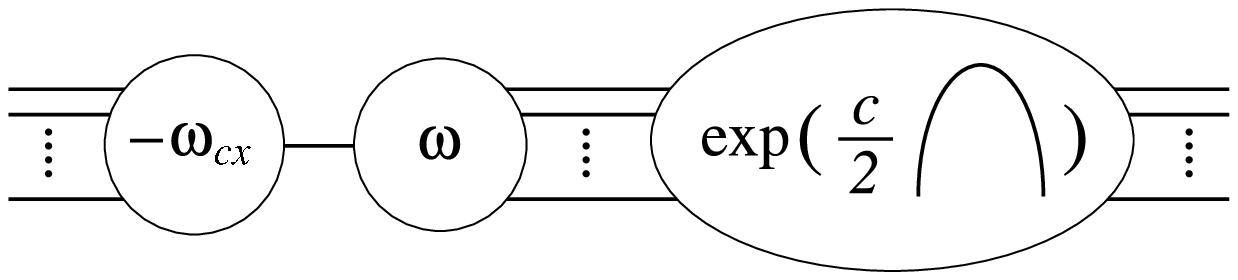}
\equiv -c\!\!
\begin{array}{c}\begin{picture}(90,30)
\put(0,5){\pc{megane3}}
\put(8,30){$f(c x)$}
\put(57,30){$f(c x)$}
\end{picture}\end{array}.
$$
This implies that $\partial_\Omega$ takes the right hand side
of the formula of the lemma to 
$\exp \big( \frac{c}2 \pc{arc} \big)$.
\end{proof}

\begin{ex}
\label{ex.p2torus}
For the $(p,2)$ torus knot,
Theorem \ref{thm.2ltorus} implies that
\begin{align*}
\Theta_{T(p,2)}(t_1,t_2,t_3) &=
\frac{1}{(t_1+1)(t_2+1)(t_3+1)} \\*
&\quad\times\Big(
\frac{p-1}2 \big( t_1^p +t_1^{-p} +t_2^p +t_2^{-p} +t_3^p +t_3^{-p} \big) \\*
&\qquad
-\frac{t_1^{p-1}-t_1^{-(p-1)}}{t_1-t_1^{-1}}
-\frac{t_2^{p-1}-t_2^{-(p-1)}}{t_2-t_2^{-1}}
-\frac{t_3^{p-1}-t_3^{-(p-1)}}{t_3-t_3^{-1}} 
\Big).
\end{align*}
For example, the coefficients of $\Theta_{T(7,2)}(t_1,t_2,t_3)$
are as shown in Table \ref{tbl.T72}.
Further,
\begin{align*}
\hat\Theta_{T(p,2)}(t) &=
\frac{t^2}{(t^2-1)^2} 
\Big( 
\frac{p-1}2 \big( t^p +t^{-p} \big) 
-\frac{t^{p-1}-t^{-(p-1)}}{t-t^{-1}}
\Big) \\
&=\frac{t^3}{(t^2-1)^3}
\Big(
\frac{p-1}2 ( t^{p+1}-t^{-p-1} ) - \frac{p+1}2 ( t^{p-1}-t^{-p+1} )
\Big).
\end{align*}
\end{ex}

\begin{table}[htpb]
$$
{\small
\begin{array}{lrrrrrrrrrrrrrr}
n &\ & -6&-5&-4&-3&-2&-1&\ \ 0&\ \ 1&2&\ \ 3&4&\ \ 5&6 \\
\\
m=6 && \cdot&\cdot&\cdot&\cdot&\cdot&\cdot& 3&-3& 3&-3& 3&-3& 3 \\
m=5 && \cdot&\cdot&\cdot&\cdot&\cdot&-3&\cdot&\cdot&\cdot&\cdot&\cdot&\cdot&-3 \\
m=4 && \cdot&\cdot&\cdot&\cdot& 3&\cdot& 2&-2& 2&-2& 2&\cdot& 3 \\
m=3 && \cdot&\cdot&\cdot&-3&\cdot&-2&\cdot&\cdot&\cdot&\cdot&-2&\cdot&-3 \\
m=2 && \cdot&\cdot& 3&\cdot& 2&\cdot& 1&-1& 1&\cdot& 2&\cdot& 3 \\
m=1 && \cdot&-3&\cdot&-2&\cdot&-1&\cdot&\cdot&-1&\cdot&-2&\cdot&-3 \\
m=0 &&  3&\cdot& 2&\cdot& 1&\cdot&\cdot&\cdot& 1&\cdot& 2&\cdot& 3 \\
m=-1 && -3&\cdot&-2&\cdot&-1&\cdot&\cdot&-1&\cdot&-2&\cdot&-3&\cdot \\
m=-2 &&  3&\cdot& 2&\cdot& 1&-1& 1&\cdot& 2&\cdot& 3&\cdot&\cdot \\
m=-3 && -3&\cdot&-2&\cdot&\cdot&\cdot&\cdot&-2&\cdot&-3&\cdot&\cdot&\cdot \\
m=-4 &&  3&\cdot& 2&-2& 2&-2& 2&\cdot& 3&\cdot&\cdot&\cdot&\cdot \\
m=-5 && -3&\cdot&\cdot&\cdot&\cdot&\cdot&\cdot&-3&\cdot&\cdot&\cdot&\cdot&\cdot \\
m=-6 &&  3&-3& 3&-3& 3&-3& 3&\cdot&\cdot&\cdot&\cdot&\cdot&\cdot
\end{array}}
$$
\caption{\label{tbl.T72}
The non-zero coefficients of $t_1^{n}t_2^{m}$ in 
$\Theta_{T(7,2)}(t_1,t_2,t_1^{-1}t_2^{-1})$}
\end{table}

\begin{proof}
By definition,
\begin{align*}
& \Delta_{T(p,2)}(t) = \frac{t^{p/2}+t^{-p/2}}{t^{1/2}+t^{-1/2}}, \quad
\psi_{p,2}(t)=-\frac{t^{p/2}-t^{-p/2}}{t^{1/2}+t^{-1/2}}, \quad \\
& \psi_{2,p}(t)=\frac{1}{(t^{1/2}+t^{-1/2})(t^{p/2}-t^{-p/2})} \cdot
\Big( (t+t^{-1})\cdot \frac{t^p-t^{-p}}{t-t^{-1}} - p(t^p+t^{-p}) \Big).
\end{align*}
Hence, when $\{ i,j,k \} = \{ 1,2,3 \}$, we have that
$$
\frac12 \Big( \psi_{p,2}(t_i) \Delta_{T(p,2)}(t_k)
+\psi_{p,2}(t_k) \Delta_{T(p,2)}(t_i) \Big)
=
\frac{t_j^{p/2}-t_j^{-p/2}}{(t_i^{1/2}+t_i^{-1/2})(t_k^{1/2}+t_k^{-1/2})} .
$$
Therefore,
\begin{align*}
&-\frac14 \psi_{2,p}(t_j) \cdot
\Big( \psi_{p,2}(t_i) \Delta_{T(p,2)}(t_k)
+\psi_{p,2}(t_k) \Delta_{T(p,2)}(t_i) \Big) \\
&=
\frac{1}{(t_i^{1/2}+t_i^{-1/2})(t_j^{1/2}+t_j^{-1/2})(t_k^{1/2}+t_k^{-1/2})} 
\cdot \frac12 \cdot
\Big( 
p(t_j^p+t_j^{-p}) 
-(t_j+t_j^{-1})\cdot \frac{t_j^p-t_j^{-p}}{t_j-t_j^{-1}} 
\Big) \\
&=
\frac{1}{(t_i^{1/2}+t_i^{-1/2})(t_j^{1/2}+t_j^{-1/2})(t_k^{1/2}+t_k^{-1/2})} 
\cdot \Big( 
\frac{p-1}2 \big( t_j^p +t_j^{-p} \big) 
-\frac{t_j^{p-1}-t_j^{-(p-1)}}{t_j-t_j^{-1}}
\Big).
\end{align*}
By Theorem \ref{thm.2ltorus},
we obtain 
$\Theta_{T(p,2)}(t_1,t_2,t_3)$
as the sum of the above formula
over $(i,j,k)=(1,2,3), (2,3,1), (3,1,2)$,
which gives the required formula.
\end{proof}

\begin{ex}
In a similar way as the previous example,
we have that
\begin{align*}
\Theta_{T(p,3)}(t_1,t_2,t_3) &=
\frac{(t_1-1)(t_2-1)(t_3-1)}{(t_1^3-1)(t_2^3-1)(t_3^3-1)} \\*
&\quad\times\Big(
(p-1) \big(
t_1^p +t_1^{-p} +t_2^p +t_2^{-p} +t_3^p +t_3^{-p} \\*
&\qquad\qquad\qquad
+t_1^{2p} +t_1^{-2p} +t_2^{2p} +t_2^{-2p} +t_3^{2p} +t_3^{-2p} \\*
&\qquad\qquad\qquad
+t_1^{2p}t_2^p +t_1^{-2p}t_2^{-p} 
+t_1^{p}t_2^{2p} +t_1^{-p}t_2^{-2p} 
+t_1^{p}t_2^{-p} +t_1^{-p}t_2^{p} \big) \\*
&\qquad-
\frac{t_1^{3(p-1)/2}-t_1^{-3(p-1)/2}}{t_1^{3/2}-t_1^{-3/2}} \cdot
\big(2t_1^{p/2} +2t_1^{-p/2} +t_2^{p/2}t_3^{-p/2} +t_2^{-p/2}t_3^{p/2} \big)\\*
&\qquad-
\frac{t_2^{3(p-1)/2}-t_2^{-3(p-1)/2}}{t_2^{3/2}-t_2^{-3/2}} \cdot
\big(2t_2^{p/2} +2t_2^{-p/2} +t_1^{p/2}t_3^{-p/2} +t_1^{-p/2}t_3^{p/2} \big)\\*
&\qquad-
\frac{t_3^{3(p-1)/2}-t_3^{-3(p-1)/2}}{t_3^{3/2}-t_3^{-3/2}} \cdot
\big(2t_3^{p/2} +2t_3^{-p/2} +t_1^{p/2}t_2^{-p/2} +t_1^{-p/2}t_2^{p/2} \big)
\Big), 
\end{align*}
and
\begin{align*}
&\hat\Theta_{T(p,3)}(t) 
= \frac{t^3(t^{p/2}+t^{-p/2})}{(t^3-1)^2} \cdot
\Big( (p-1)(t^{3p/2}+t^{-3p/2})
- 2 \cdot \frac{t^{3(p-1)/2}-t^{-3(p-1)/2}}{t^{3/2}-t^{-3/2}} 
\Big) \\*
&= \frac{t^{p/2}+t^{-p/2}}{(t^{3/2}-t^{-3/2})^3} \cdot
\Big( (p-1)( t^{3(p+1)/2} -t^{-3(p+1)/2} )
-(p+1)( t^{3(p-1)/2} -t^{-3(p-1)/2} ) \Big).
\end{align*}
See also Tables \ref{tbl.theta_pq} and \ref{tbl.hT_pq}
for the values of $\Theta_{T(p,q)}$ and $\hat\Theta_{T(p,q)}$
for some $(p,q)$.
\end{ex}

\begin{table}[htpb]
$$
\begin{array}{rl}
(p,q): &
\mbox{\small The non-zero coefficients of $t_1^n t_2^m$ in 
$\Theta_{T(p,q)}(t_1,t_2,t_1^{-1}t_2^{-1})$
in the fundamental domain}
\\
\\
(3,2): &
{\small
\begin{array}{rrr}
&&-1 \\
\ \cdot&\ \cdot& 1
\end{array}} 
\\
(5,2): &
{\small
\begin{array}{rrrrr}
&&&& 2 \\
&&-1&\ \ \ \cdot&-2 \\
\ \cdot&\ \cdot& 1&\ \cdot& 2
\end{array}}
\\
(7,2): &
{\small
\begin{array}{rrrrrrr}
&&&&&&-3 \\
&&&& 2&\cdot& 3 \\
&&-1&\ \ \ \cdot&-2&\ \ \ \cdot&-3 \\
\ \cdot&\ \cdot& 1&\ \cdot& 2&\ \cdot& 3 
\end{array}}
\\
\\
(4,3): &
{\small
\begin{array}{rrrrrrr}
&&&&&& 3 \\
&&&& 3&-3&\cdot \\
&& 1&-2&\cdot& 3&-3 \\
\ \cdot&\ \cdot&-1& 4&-3&\cdot& 3 
\end{array}}
\\
(5,3): &
{\small
\begin{array}{rrrrrrrrr}
&&&&&&&&-4 \\
&&&&&&-4&\cdot& 4 \\
&&&&-6& 3& 4&-4&\cdot \\
&&-2& 1& 3&-6&\cdot& 4&-4 \\
\ \cdot&\ \cdot& 2&-2&\cdot& 6&-4&\cdot& 4 \\
\end{array}}
\\
(7,3): &
{\small
\begin{array}{rrrrrrrrrrrrr}
&&&&&&&&&&&& 6 \\
&&&&&&&&&& 6&-6&\cdot \\
&&&&&&&& 10&-5&\cdot& 6&-6 \\
&&&&&& 12&-5&-5& 10&-6&\cdot& 6 \\
&&&& 6&-5&-4& 10&\!\!\!-10&\cdot& 6&-6&\cdot \\
&& 2&-3&-1& 6&-8&\cdot& 10&\!\!\!-10&\cdot& 6&-6 \\
\ \cdot&\ \cdot&-2& 6&-4&-2& 12&\!\!\!-10&\cdot& 10&-6&\cdot& 6 
\end{array}}
\\
\\
(5,4): &
{\small
\begin{array}{rrrrrrrrrrrrr}
&&&&&&&&&&&&-6 \\
&&&&&&&&&&-6& 6&\cdot \\
&&&&&&&& 9&\cdot&\cdot&-6& 6 \\
&&&&&& 1&-5&\cdot&\cdot& 6&-6&\cdot \\
&&&&-5& 4&-4& 5&-5&\cdot&\cdot&\cdot&\cdot \\
&& 1& 1&-2&\cdot& 3& 1&-4&\ \ \ \cdot&\cdot& 6&-6 \\
\ \cdot&\ \cdot&-1&-2& 9&-8& 1&-2& 9&\cdot&-6&\cdot& 6 
\end{array}}
\\
(7,4): &
{\small
\begin{array}{rrrrrrrrrrrrrrrrrrr}
&&&&&&&&&&&&&&&&&&-9 \\
&&&&&&&&&&&&&&&&-9& 9&\cdot \\
&&&&&&&&&&&&&& 15&\ \ \cdot&\cdot&-9& 9 \\
&&&&&&&&&&&& 15&\!\!\!-15&\cdot&\cdot& 9&-9&\cdot \\
&&&&&&&&&&-6& -6&\cdot& 15&\!\!\!-15&\cdot&\cdot&9&-9 \\
&&&&&&&&\!\!\!-18& 7&-1&12&\!\!\!-15& 8&7&\cdot&-9& \cdot&9 \\
&&&&&&-8&10&5&-6& 1&\cdot& 7&\!\!\!-15& 8&\cdot&9&-9& \cdot \\
&&&&5&-5&4&\!\!\!-11&13&\!\!\!-13& 7&-7& 8&-8& \cdot&\cdot&\cdot&\cdot&\cdot \\
&&2&-4&2&4&2&-9&\cdot&11&-6&-5&\cdot&15&\!\!\!-15&\cdot&\cdot&9&-9 \\
\ \cdot&\cdot&-2&8&-9&2&-4& 20&\!\!\!-18&2&4&12&\!\!\!-15&\cdot&15&\cdot&-9&\cdot&9 
\end{array}}
\end{array}
$$
\caption{\label{tbl.theta_pq}
The non-zero coefficients of $t_1^n t_2^m$ in 
$\Theta_{T(p,q)}(t_1,t_2,t_1^{-1}t_2^{-1})$
in a fundamental domain $\{ 0 \le 2m \le n \}$ (see \cite{Rozansky_2lp})
for $(p,q)$ with $p\le7$, $q\le4$.
The array for each $(p,q)$
is a subset of the full array such as shown in Table \ref{tbl.T72}
and the most left dot is at $(n,m)=(0,0)$.
We can recover the other coefficients for each $(p,q)$
from the presented coefficients by the symmetry of 
$\Theta_K(t_1,t_2,t_1^{-1}t_2^{-1})$.}
\end{table}

\begin{table}[htpb]
$$
\small
\begin{array}{ll}
(p,q): & \mbox{The part of non-negative powers in $\hat\Theta_{T(p,q)}(t)$} \\
\\
(3,2): &
t
\\
(5,2): &
3 t + 2 t^3
\\
(7,2): &
6 t + 5 t^3  + 3 t^5
\\
(9,2): &
10 t + 9 t^3  + 7 t^5  + 4 t^7
\\
\\
(4,3): &
3 t + 4 t^2  + 3 t^5
\\
(5,3): &
6 t + 4 t^2  + 6 t^4  + 4 t^7
\\
(7,3): &
10 t + 12 t^2  + 6 t^4  + 12 t^5  + 10 t^8  + 6 t^{11}
\\
(8,3): &
15 t + 12 t^2  + 16 t^4  + 7 t^5  + 15 t^7  + 12 t^{10} + 7 t^{13}
\\
(10,3): &
21 t +24 t^2 +16 t^4 +25 t^5 +9 t^7 +24 t^8 +21 t^{11} +16 t^{14} + 9t^{17}
\\
\\
(5,4): &
6 t + 12 t^2  + 9 t^3  + 8 t^6  + 9 t^7  + 6 t^{11}
\\
(7,4): &
15t +24t^2 + 9t^3 +18t^5 +20t^6 +18t^9 +12t^{10} +15t^{13} +9t^{17}
\\
(9,4): &
21t +40t^2 + 27t^3 +12t^5 +36t^6 +30t^7 +28t^{10} +30t^{11} +16t^{14}   
+27t^{15} +21t^{19} +12t^{23}
\\
\\
(6,5): &
10t +24t^2 + 27t^3 +16t^4 +15t^7 +24t^8 +18t^9 +15t^{13} +16t^{14} +10t^{19}
\\
(7,5): &
36t +12t^2 +20t^3 +30t^4 +36t^6 +24t^8 +18t^9 +30t^{11} +24t^{13} 
+18t^{16} +20t^{18} +12t^{23}
\\
(8,5): &
45 t + 24 t^2  + 14 t^3  + 48 t^4  + 36 t^6  + 30 t^7  + 45 t^9  
+ 21t^{11} +32t^{12} +36t^{14} +30t^{17} +21t^{19} +24t^{22} +14t^{27}
\\
(9,5): &
28 t \!+\! 60 t^2  \!+\! 54 t^3  \!+\! 16 t^4  \!+\! 36 t^6  \!+\! 60 t^7  
\!+\!42 t^8  \!+\! 40 t^{11} \!\!+\! 54 t^{12} \!\!+\! 24 t^{13} \!\!+\! 40 t^{16}   
\!\!+\! 42 t^{17} \!\!+\! 36 t^{21} \!\!+\! 24 t^{22} \!\!+\!   28 t^{26} \!\!+\! 16 t^{31}
\end{array}
$$
\caption{\label{tbl.hT_pq}
The parts of non-negative powers in $\hat\Theta_{T(p,q)}(t)$
for $(p,q)$ with $p\le10$, $q\le5$.
The remaining part for each $(p,q)$ can recover from the presented part
by replacing $t$ with $t^{-1}$.}
\end{table}

\section{A cabling formula for the 2-loop polynomial}

In this section,
we give a cabling formula for the 2-loop polynomial.
We show the formula by picking up the 2-loop part of 
the cabling formula of the Kontsevich invariant,
modifying the proof of Theorem \ref{thm.2ltorus}.
This cabling formula is also obtained independently by 
March\'e \cite{Marche_PhD}.

It is known, see {\it e.g.} \cite{Lickorish_book},
that a cabling formula for the Alexander polynomial is given by
$$
\Delta_{K^{(p,q)}}(t) = \Delta_{T(p,q)}(t) \Delta_K(t^p).
$$
A cabling formula for the 2-loop polynomial is given by

\begin{thm}
\label{thm.cabling}
Let $K$ be a knot,
and let $K^{(p,q)}$ be the $(p,q)$ cable knot of $K$.
Then, 
\begin{align*}
&\Theta_{K^{(p,q)}}(t_1,t_2,t_3) =
\Theta_{T(p,q)}(t_1,t_2,t_3)
+ \Theta_K(t_1^p,t_2^p,t_3^p) \\
&\quad+ \frac12
\Delta_{T(p,q)}(t_1) \Delta_{T(p,q)}(t_2) \Delta_{T(p,q)}(t_3)
\!\!\!\!\!\! \sum_{ \{i,j,k\}=\{1,2,3\} } \!\!\!\!\!\!
\Delta'_K(t_i^p) \cdot t_i^p 
\cdot \phi_{q,p}(t_j) \Delta_K(t_j^p) \Delta_{K}(t_k^p).
\end{align*}
\end{thm}

\begin{proof}
We show the theorem, modifying the proof of Theorem \ref{thm.2ltorus}.
By Proposition \ref{prop.cabling}, we have that
$$
Z^w\big( K^{(p,q)} \big)
\equiv
\partial_\Omega^{-1} \Psi^{(p)} \partial_\Omega
\Big( Z^w(K) \sqcup
\exp_\sqcup \big( \frac{q}{2p} \pc{arc} \big) \Big)
\sqcup 
\exp_\sqcup \big( -\frac{p q}{2} \pc{arc} \big) ,
$$
where $Z^w(K)$ is presented by
$$
Z^w(K) = \Omega \sqcup
\exp_\sqcup \Big( 
{\begin{array}{c}
{\begin{picture}(62,10)
\put(-5,20){$-\frac12 \log \Delta{\!}_K{\!}(e^x)$}
\put(-3,-5){\pc{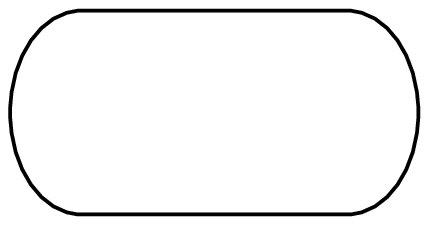}}
\end{picture}} \end{array}}  \Big)
+ (\mbox{terms of $(\ge2)$-loop}).
$$
The 2-loop part of $Z^w(K)$ contributes to the required formula by
$\Theta_K(t_1^p,t_2^p,t_3^p)$.
We calculate the contribution of the 1-loop part
in the following of this proof.

In a similar way as (\ref{eq.dOeO}), we have that
\begin{align*}
& \partial_\Omega \Big( Z^w(K) \sqcup
\exp_\sqcup \big( \frac{q}{2p} \pc{arc} \big) 
\Big) \\*
&\equiv
\exp \big( \frac{q}{2p} \pc{arc} \big)
\sqcup \Omega \sqcup \Omega_{\frac{q}{p} x} \sqcup 
\exp_\sqcup \Big( 
{\begin{array}{c}
{\begin{picture}(62,10)
\put(-5,20){$-\frac12 \log \Delta{\!}_K{\!}(e^x)$}
\put(-3,-5){\pc{loop1}}
\end{picture}} \end{array}}
+
\! \begin{array}{c}\begin{picture}(100,10)
\put(0,-5){\pc{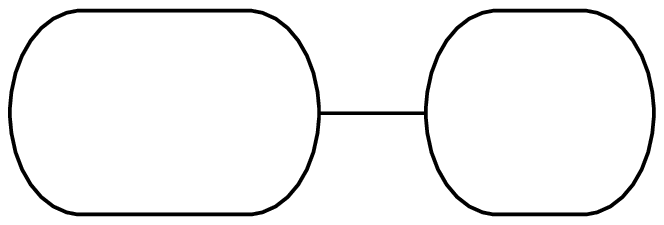}}
\put(0,20){$f(x)\!+\!g(x)$}
\put(65,22){$f( \frac{q}{p} x)$}
\end{picture}\end{array} \!
\Big),
\end{align*}
where $g(x)$ is given by
$$
g(x) = \frac{d}{dx} \Big( -\frac12 \log \Delta_K(e^x) \Big)
= -\frac{\Delta'_K(e^x) \cdot e^x}{2 \Delta_K(e^x)}.
$$
The map $\Psi^{(p)}$ sends this to
$$
\exp \big( \frac{p q}{2} \pc{arc} \big)
\sqcup \Omega_{p x} \sqcup \Omega_{q x} \sqcup 
\exp_\sqcup \Big( 
{\begin{array}{c}
{\begin{picture}(62,10)
\put(-5,20){$-\frac12 \log \Delta{\!}_K{\!}(e^{p x})$}
\put(-3,-5){\pc{loop1}}
\end{picture}} \end{array}}
+
\! \begin{array}{c}\begin{picture}(100,10)
\put(0,-5){\pc{megane1}}
\put(-10,20){$f(p x)\!+\!g(p x)$}
\put(65,20){$f(q x)$}
\end{picture}\end{array} \!
\Big).
$$
Calculating its image by $\partial_\Omega^{-1}$
in a similar way as in the proof of Theorem \ref{thm.2ltorus},
the error term coresponding to the formula (\ref{eq.me-theta})
is as follows,
\begin{align*}
& \begin{array}{c}\begin{picture}(90,35)
\put(0,5){\pc{megane3}}
\put(8,30){$g(p x)$}
\put(57,30){$f(q x)$}
\end{picture}\end{array}
-p \!\!
\begin{array}{c}\begin{picture}(90,35)
\put(0,5){\pc{megane3}}
\put(8,30){$g(p x)$}
\put(52,30){$f(p q x)$}
\end{picture}\end{array} \\*
&=
\frac1{4}
\begin{array}{c}
\begin{picture}(110,30)
\put(0,5){\pc{megane2}}
\put(10,32){$g(p x)$}
\put(70,32){$\phi_{q,p}(t)$}
\end{picture}\end{array}
= - \frac12
\begin{array}{c}
\begin{picture}(80,65)
\put(0,20){\pc{theta2}}
\put(22,60){$g(p x)$}
\put(22,32){$\phi_{q,p}(t)$}
\end{picture}\end{array} .
\end{align*}
This contributes to the required formula by
$$
\sum_{ \{i,j,k\}=\{1,2,3\} } \!\!\!
\frac{\Delta'_K(t_i^p) \cdot t_i^p}{2 \Delta_K(t_i^p)}
\cdot \Delta_{K^{(p,q)}}(t_i)
\phi_{q,p}(t_j) \Delta_{K^{(p,q)}}(t_j) \Delta_{K^{(p,q)}}(t_k).
$$
Noting that
$\Delta_{K^{(p,q)}}(t) = \Delta_{T(p,q)}(t) \Delta_K(t^p)$,
we obtain the required formula.
\end{proof}

A cabling formula for the reduced 2-loop polynomial is given by

\begin{cor}
\label{cor.chT}
For the notation in Theorem \ref{thm.cabling},
\begin{align*}
\hat\Theta_{K^{(p,q)}}(t)
&= 
\hat\Theta_{T(p,q)}(t)
+ \frac{(t^{p/2}-t^{-p/2})^2}{(t^{1/2}-t^{-1/2})^2} \cdot
\hat\Theta_K(t^p) \\*
&\quad -\frac{t^p}{(t^{1/2}-t^{-1/2})^2} \cdot
\Delta_{T(p,q)}(t) 
\Delta_K(t^p) \Delta'_K(t^p) \psi_{q,p}(t) .
\end{align*}
\end{cor}

\begin{proof}
The required formula is obtained from the formula of Theorem \ref{thm.cabling}
by putting $t_1=t$,\ \ $t_2=1/t$, and $t_3=1$.
\end{proof}

\section{Relations to Vassiliev invariants}

In this section
we show some relations to Vassiliev invariants of degree 2, 3.

A leading part of the Kontsevich invariant is presented by
$$
\log_\sqcup Z^\sigma(K) - \omega
\ \ =\ \ \frac{v_2(K)}2 \pc{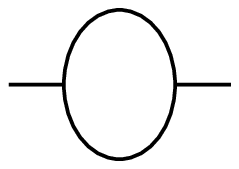}
+ \frac{v_3(K)}4 \pc{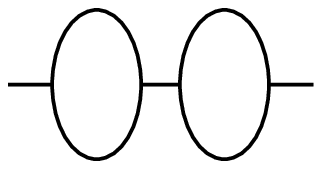}
+ (\mbox{terms of degree $\ge 4$}),
$$
where the {\it degree} of a Jacobi diagram is
half the number of univalent and trivalent vertices of the diagram,
and $v_2$, $v_3$ are $\Z$-valued primitive Vassiliev invariants of degree 2, 3
respectively (see \cite{Ohtsuki_book}).
Since $\pc{v2}$ has 1-loop,\vspace{-0.2pc}
$v_2(K)$ can be presented by the Alexander polynomial;
in fact, from the formula of the loop expansion, 
\begin{align*}
v_2(K) &= - \big(
\mbox{the coefficient of $x^2$ in the expansion of $\Delta_K(e^x)$} \big) \\*
&= -\frac12 \Delta''_K(1).
\end{align*}
Further, since $\pc{v3}$ has 2-loop,\vspace{-0.3pc}
$v_3(K)$ can be presented by the 2-loop polynomial; in fact, we have

\begin{prop}
\label{prop.v3_Theta1}
$$
v_3(K) = \frac12 \hat\Theta_K(1) .
$$
\end{prop}

\begin{proof}
Let us consider the map
\begin{align*}
\pc{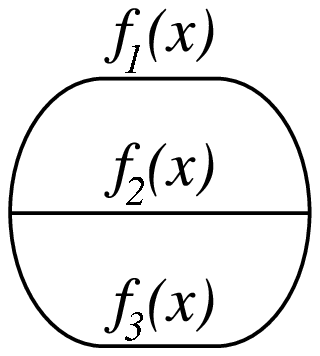} \ &\longmapsto \ 
f_3(0)
\begin{picture}(50,20)
\put(0,10){\pc{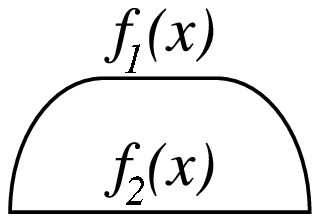}}
\end{picture}
\ +\ f_2(0)
\pc{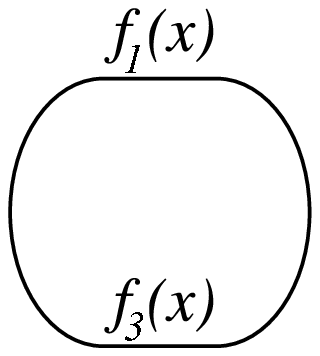} 
\ +\ f_1(0)
\begin{picture}(50,20)
\put(0,-9){\pc{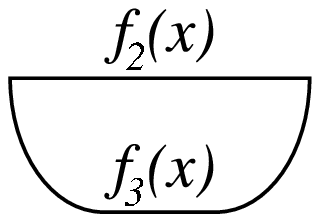}}
\end{picture} \\
&\longmapsto
\frac16 \!\sum_{ \{i,j,k\}=\{1,2,3\} }\!\!\!
f_i(x) f_j(-x) f_k(0).
\end{align*}
This map takes the 2-loop part of $\log_\sqcup Z^\sigma(K)$
to $\frac1{12}(e^{x/2}-e^{-x/2})^2 \hat\Theta_K(e^x)$,
whose coefficient of $x^2$ equals $\frac1{12} \hat\Theta_K(1)$.
Since $\pc{v3}=\pc{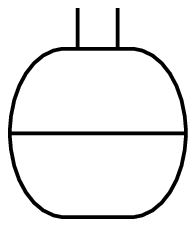}$ by the AS and IHX relations,
the above maps takes this diagram to $\frac23 x^2$.
Hence, $\frac16 v_3(K) = \frac1{12}\hat\Theta_K(1)$,
which implies the required formula.
\end{proof}

\begin{ex}
A cabling formula for $v_3$ is given by
$$
v_3(K^{(p,q)}) = p^2 \cdot v_3(K) 
+ \frac1{12} p(p^2-1)q \cdot \Delta''_K(1) 
+ \frac1{144} p(p^2-1)q(q^2-1).
$$
\end{ex}

\begin{proof}
From Proposition \ref{prop.v3_Theta1}
and Corollary \ref{cor.chT} putting $t=1$,
we have that
$$
v_3(K^{(p,q)}) = v_3\big( T(p,q) \big) + p^2 \cdot v_3(K) 
- \frac{p}2 \Delta''_K(1) \phi'_{q,p}(1).
$$
The required formula follows from it, by using
\begin{align*}
& v_3 \big( T(p,q) \big)
= \frac12 \hat\Theta_{T(p,q)}(1) = \frac1{144} p(p^2-1)q(q^2-1), \\
& \phi'_{q,p}(1) = \frac16 q (1-p^2).
\end{align*}
For the value of the first formula, see also \cite{Willerton}.
\end{proof}

\normalsize

Research Institute for Mathematical Sciences, 
Kyoto University, 
Sakyo-ku, Kyoto, 606-8502, 
Japan

E-mail address: tomotada@kurims.kyoto-u.ac.jp

\end{document}